# Continuity of Metric Projection Operator from $C[0, 1]$ onto $\mathcal{P}_n$ with Applications to Mordukhovich Derivatives

Jinlu Li

Department of Mathematics
Shawnee State University
Portsmouth, Ohio 45662 USA
jli@shawnee.edu

**Abstract**

Let $C[0, 1]$ be the Banach space of all continuous real valued functions on [0, 1]. For an arbitrarily given nonnegative integer $n$, let $\mathcal{P}_n$ denote the set of all polynomials with degree less than or equal to $n$. $\mathcal{P}_n$ is a closed subspace of $C[0, 1]$. In this paper, we first prove (in details) that the metric projection $P_{\mathcal{P}_n}$: $C[0, 1]$ $\rightarrow \mathcal{P}_n$ is a single-valued mapping and it is (norm to norm) continuous. Then, we use the continuity of $\mathcal{P}_n$ to investigate the Gâteaux directional derivatives, some properties and fixed-point properties of the Mordukhovich derivatives of $\mathcal{P}_n$.

## 1. Introduction and Preliminaries

Let $(X, \|\cdot\|)$ be a Banach space and $C$ a nonempty closed and convex subset of $X$. Let $P_C$ denote the (standard) metric projection operator from $X$ onto $C$, which is defined, for any $x \in X$, by

$$P_C x = \{y \in C: \|x - y\| \leq \|x - z\|, \text{for all } z \in C\}.$$

**Remarks 1.1** (see [1, 2, 9, 12]). It is well-known that if $X$ is a uniformly convex and uniformly smooth Banach space, then $P_C: X \rightarrow C$ is a single-valued continuous mapping. In general, $P_C: X \rightarrow 2^C\backslash\{\emptyset\}$ is indeed a set-valued mapping.

The normalized duality mapping $J: X \rightrightarrows X^*$ is a set-valued mapping such that,

$$J(x) = \{\varphi \in X^*: \langle\varphi, x\rangle = \|x\|^2 = \|\varphi\|_*^2\}, \text{ for any } x \in X.$$

For the properties of the normalized duality mapping, one may see Sections 4.2–4.3 and Problem set 4.2 in [28] and [1, 2, 9, 12, 17, 23, 24, 25, 30].

Let $rca[0, 1]$ be the Banach space of all real valued, regular and countable additive functionals on [0, 1], in which the considered σ-field is the standard Borel σ-field Σ on [0, 1] including all closed and open subsets of [0, 1]. $rca[0, 1]$ is equipped with the norm of total variation. That is, for $\mu \in rca[0, 1]$, the norm of $\mu$ is defined by its total variation $v(\mu, [0, 1])$ (see page 160 in [7]).

Let ($C$[0, 1], ‖·‖) be the Banach space of all continuous real valued functions on [0, 1] with respect to the standard Borel $\sigma$-field Σ and with the maximum norm

$$\|f\| = \max_{0\le t\le 1}|f(t)|, \text{ for any } f \in C[0, 1].$$

The dual space of $C$[0, 1] is denoted by $C^*[0, 1]$ with norm $\|\cdot\|_*$. Let $\langle\cdot,\ \cdot\rangle$ denote the real canonical pairing between $C^*[0,1]$ and $C$[0, 1]. By the Riesz Representation Theorem (see Theorem IV.6.3 in [7]), for any φ ∈ $C^*[0,1]$ , there is μ ∈ $rca$[0, 1], such that

$$\langle \varphi, f\rangle = \int_0^1 f(t)\mu(dt), \text{ for any } f \in C[0, 1]. \tag{1.1}$$

Then, throughout this paper, we say that $C^*[0, 1] = rca[0, 1]$. Without any special mention, we shall identify φ and $\mu$ in (1.1), and (1.1) is rewritten as

$$\langle \mu, f\rangle = \int_0^1 f(t)\mu(dt), \text{ for any } f \in C[0, 1].$$

It satisfies

$$\|\mu\|_* = v(\mu, [0, 1]) = \sup_{E\in\Sigma}|\mu(E)|.$$

The origin of $C$[0, 1] is denoted by $\theta$, which is the constant function defined on [0, 1] with value 0. The origin of the dual space $C^*[0, 1]$ is denoted by $\theta^*$, which is also a constant functional on Σ with value 0. This is, $\theta^*(E) = 0$, for every $E \in \Sigma$. This implies

$$\|\theta^*\|_* = v(\theta^*, [0, 1]) = \sup_{E\in\Sigma}|\theta^*(E)| = 0.$$

For any given nonnegative integer $n$, let $\mathcal{P}_n$ denote the closed subspace of $C$[0, 1] that consists of all real coefficient polynomials of degree less than or equal to $n$. $P_{\mathcal{P}_n}$ is the metric projection from $C$[0, 1] to $\mathcal{P}_n$. It is well-known that $C$[0, 1] is a nonreflexive Banach space. Therefore $C$[0, 1] is not a uniformly convex and uniformly smooth Banach space, which follows that the single-valued property and the continuity of $P_{\mathcal{P}_n}$ cannot be automatically implied by Remarks 1.1. However, in section 2 of this paper, we prove some properties of bounded subsets in $\mathcal{P}_n$. By using these properties and by using the Chebyshev's Equioscillation Theorem, we provide a new proof for $P_{\mathcal{P}_n}$ to be a single-valued mapping in section 3. In section 4, we prove that $P_{\mathcal{P}_n}$ is a single-valued continuous mapping. As applications of the continuity of $P_{\mathcal{P}_n}$, in section 5, we investigate the Gâteaux directional derivatives of $\mathcal{P}_n$ and some basic properties including some fixed-point properties of the Mordukhovich derivatives $\mathcal{P}_n$.

## 2. Properties of Bounded Subsets in $\mathcal{P}_n$

In this section, we prove some properties of bounded subsets in $\mathcal{P}_n$. More precisely speaking, we prove some properties of the boundness of the coefficients of polynomials in a bounded subset of $\mathcal{P}_n$, which will be used in the proof of the single-valued property in section 3, and the continuity of $P_{\mathcal{P}_n}$ in section 4.

**Lemma 2.1**. *For any positive integer n and for* $t \in [0, 1]$, *define* $A_n \in C[0, 1]$ *by the following* $(n+1)\times(n+1)$ *determinant*

$$A_n(t) = \begin{vmatrix} 1 & 1 & 1 & 1 & \dots & 1 & 1 \\ 1 & t & t^2 & t^3 & \dots & t^{n-1} & t^n \\ 1 & t^2 & t^4 & t^6 & \dots & t^{2(n-1)} & t^{2n} \\ 1 & t^3 & t^6 & t^9 & \dots & t^{3(n-1)} & t^{3n} \\ \dots & \dots & \dots & \dots & \dots & \dots & \dots \\ 1 & t^{n-1} & t^{2(n-1)} & t^{3(n-1)} & \dots & t^{(n-1)^2} & t^{(n-1)n} \\ 1 & t^n & t^{2n} & t^{3n} & \dots & t^{(n-1)n} & t^{n^2} \end{vmatrix}.$$

*Then,*

$$A_n(t) \neq 0, \text{ for every } t \in (0, 1).$$

*Proof*. We proceed by induction.

Step 1. For $n = 2$. We have

$$A_2(t) = \begin{vmatrix} 1 & 1 & 1 \\ 1 & t & t^2 \\ 1 & t^2 & t^4 \end{vmatrix}.$$

One can check that

$$A_2(t) = t^5 - 2t^4 + 2t^2 - t = t(t^4 - 2t^3 + 2t - 1) < 0, \text{ for every } t \in (0, 1).$$

Hence, the lemma is true for $n = 2$.

Step 2. Assume that, for $n > 2$, we have

$$A_{n-1}(t) \neq 0, \text{ for every } t \in (0, 1).$$

Step 3. We prove $A_n(t) \neq 0$, for every $t \in (0, 1)$ for $n$. For any $t \in (0, 1)$, by the properties of determinates, we have

$$A_n(t) = \begin{vmatrix} 1 & 1 & 1 & 1 & \dots & 1 & 1 \\ 1 & t & t^2 & t^3 & \dots & t^{n-1} & t^n \\ 1 & t^2 & t^4 & t^6 & \dots & t^{2(n-1)} & t^{2n} \\ 1 & t^3 & t^6 & t^9 & \dots & t^{3(n-1)} & t^{3n} \\ \dots & \dots & \dots & \dots & \dots & \dots & \dots \\ 1 & t^{n-1} & t^{2(n-1)} & t^{3(n-1)} & \dots & t^{(n-1)^2} & t^{(n-1)n} \\ 1 & t^n & t^{2n} & t^{3n} & \dots & t^{(n-1)n} & t^{n^2} \end{vmatrix}$$

$$= \begin{vmatrix} 1 & 1 & 1 & 1 & \dots & 1 & 1 \\ 0 & t-1 & t^2-1 & t^3-1 & \dots & t^{n-1}-1 & t^n-1 \\ 0 & t^2-1 & t^4-1 & t^6-1 & \dots & t^{2(n-1)}-1 & t^{2n}-1 \\ 0 & t^3-1 & t^6-1 & t^9-1 & \dots & t^{3(n-1)}-1 & t^{3n}-1 \\ \dots & \dots & \dots & \dots & \dots & \dots & \dots \\ 0 & t^{n-1}-1 & t^{2(n-1)}-1 & t^{3(n-1)}-1 & \dots & t^{(n-1)^2}-1 & t^{(n-1)n}-1 \\ 0 & t^n-1 & t^{2n}-1 & t^{3n}-1 & \dots & t^{n(n-1)}-1 & t^{n^2}-1 \end{vmatrix}$$

$$= \begin{vmatrix} 1 & 0 & 0 & 0 & \dots & 0 & 0 \\ 0 & t-1 & t^2-1 & t^3-1 & \dots & t^{n-1}-1 & t^n-1 \\ 0 & t^2-1 & t^4-1 & t^6-1 & \dots & t^{2(n-1)}-1 & t^{2n}-1 \\ 0 & t^3-1 & t^6-1 & t^9-1 & \dots & t^{3(n-1)}-1 & t^{3n}-1 \\ \dots & \dots & \dots & \dots & \dots & \dots & \dots \\ 0 & t^{n-1}-1 & t^{2(n-1)}-1 & t^{3(n-1)}-1 & \dots & t^{(n-1)^2}-1 & t^{(n-1)n}-1 \\ 0 & t^n-1 & t^{2n}-1 & t^{3n}-1 & \dots & t^{n(n-1)}-1 & t^{n^2}-1 \end{vmatrix}$$

$$= \begin{vmatrix} 1 & 0 & 0 & 0 & \dots & 0 & 0 \\ 0 & t-1 & t^2-1 & t^3-1 & \dots & t^{n-1}-1 & t^n-1 \\ 0 & t^2-t & t^4-t^2 & t^6-t^3 & \dots & t^{2(n-1)}-t^{n-1} & t^{2n}-t^n \\ 0 & t^3-t^2 & t^6-t^4 & t^9-t^6 & \dots & t^{3(n-1)}-t^{2(n-1)} & t^{3n}-t^{2n} \\ \dots & \dots & \dots & \dots & \dots & \dots & \dots \\ 0 & t^{n-1}-t^{n-2} & t^{2(n-1)}-t^{2(n-2)} & t^{3(n-1)}-t^{3(n-2)} & \dots & t^{(n-1)^2}-t^{(n-1)(n-2)} & t^{(n-1)n}-t^{(n-2)n} \\ 0 & t^n-t^{n-1} & t^{2n}-t^{2(n-1)} & t^{3n}-t^{3(n-1)} & \dots & t^{n(n-1)}-t^{(n-1)^2} & t^{n^2}-t^{(n-1)n} \end{vmatrix}$$

$$= (t-1)(t^2-1)\dots(t^{n-1}-1)(t^n-1) \begin{vmatrix} 1 & 0 & 0 & 0 & \dots & 0 & 0 \\ 0 & 1 & 1 & 1 & \dots & 1 & 1 \\ 0 & t & t^2 & t^3 & \dots & t^{n-1} & t^n \\ 0 & t^2 & t^4 & t^6 & \dots & t^{2(n-1)} & t^{2n} \\ \dots & \dots & \dots & \dots & \dots & \dots & \dots \\ 0 & t^{n-2} & t^{2(n-2)} & t^{3(n-2)} & \dots & t^{(n-2)(n-1)} & t^{(n-2)n} \\ 0 & t^{n-1} & t^{2(n-1)} & t^{3(n-1)} & \dots & t^{(n-1)^2} & t^{(n-1)n} \end{vmatrix}$$

$$= (t-1)(t^2-1)\dots(t^{n-1}-1)(t^n-1)\, t\, t^2 \dots t^{n-2} t^{n-1}$$

$$\times \begin{vmatrix} 1 & 0 & 0 & 0 & \dots & 0 & 0 \\ 0 & 1 & 1 & 1 & \dots & 1 & 1 \\ 0 & 1 & t & t^2 & \dots & t^{n-2} & t^{n-1} \\ 0 & 1 & t^2 & t^4 & \dots & t^{2(n-2)} & t^{2(n-1)} \\ \dots & \dots & \dots & \dots & \dots & \dots & \dots \\ 0 & 1 & t^{n-2} & t^{2(n-2)} & \dots & t^{(n-2)^2} & t^{(n-2)(n-1)} \\ 0 & 1 & t^{n-1} & t^{2(n-1)} & \dots & t^{(n-2)(n-1)} & t^{(n-1)^2} \end{vmatrix}$$

$$= (t-1)(t^2-1)(t^3-1)\dots(t^{n-1}-1)(t^n-1)\, t\, t^2 \dots t^{n-2} t^{n-1} A_{n-1}(t).$$

By the assumption in step 2, this implies

$$A_n(t) \neq 0, \text{for every } t \in (0, 1).$$

This proves the claim. □

Let $D$ be a nonempty subset of $\mathcal{P}_n$. Let Coe($D$) be the collection of all coefficients of polynomials in $D$, which is called the coefficient set of $D$. That is,

$$\text{Coe}(D) = \{a \in \mathbb{R}\text{: there is } p \in D \text{ such that } a \text{ is a coefficient of } p\}.$$

**Proposition 2.2**. *Let n be a positive integer. Let D be a nonempty subset of $\mathcal{P}_n$ with coefficient set* Coe($D$). *Then, we have*

$$D \textit{ is } \|\cdot\|\textit{-bounded} \quad \Leftrightarrow \quad \mathrm{Coe}(D) \textit{ is } |\cdot|\textit{-bounded}.$$

*Proof*. "⟸". Suppose that Coe($D$) is $|\cdot|$-bounded by a positive number $d$. Since we only consider the polynomials on [0, 1], this implies that $D$ is $\|\cdot\|$-bounded by $(n+1)d$.

"⟹". Suppose that $D$ is $\|\cdot\|$-bounded by $b$, for some positive number $b$. Take an arbitrary element $p \in D$ with

$$p(t) = a_0 + a_1t + a_2t^2 + \ldots + a_nt^n, \text{ for } t \in [0, 1]. \tag{2.1}$$

Respectively substituting $t$ by $1, \frac{1}{2}, \frac{1}{2^2}, \ldots, \frac{1}{2^n}$ in (2.1), we calculate the corresponding values of $p$ and we obtain the following system of linear equations with respect to $a_0, a_1, a_2, \ldots, a_n$.

$$\begin{aligned}
&a_0 + a_1 + a_2 + \ldots + a_n = b_0,\\
&a_0 + a_1\frac{1}{2} + a_2\frac{1}{2^2} + \ldots + a_n\frac{1}{2^n} = b_1,\\
&a_0 + a_1\frac{1}{2^2} + a_2\frac{1}{2^4} + \ldots + a_n\frac{1}{2^{2n}} = b_2,\\
&\ldots\ldots\\
&a_0 + a_1\frac{1}{2^n} + a_2\frac{1}{2^{2n}} + \ldots + a_n\frac{1}{2^{n^2}} = b_n.
\end{aligned} \tag{2.2}$$

Where $b_0, b_1, b_2, \ldots, b_n$ are real numbers. By the assumption that $D$ is $\|\cdot\|$-bounded by $b$, we have

$$|b_k| \le b, \text{ for } k = 0, 1, 2, \ldots, n.$$

By the definition of $A_n(\cdot)$, the system of linear equations (2.2) in terms of the "variables" $a_0, a_1, a_2, \ldots, a_n$ can be rewritten as

$$A_n\left(\frac{1}{2}\right)\begin{pmatrix} a_0 \\ a_1 \\ a_2 \\ \cdots \\ a_n \end{pmatrix} = \begin{pmatrix} b_0 \\ b_1 \\ b_2 \\ \cdots \\ b_n \end{pmatrix}.$$

For $k = 0, 1, 2, \ldots, n$, let $A_n^k\left(\frac{1}{2}\right)$ denote the $(n+1)\times(n+1)$ determinate that is defined by substituting the $k$th-column in $A_n\left(\frac{1}{2}\right)$ by $(b_0, b_1, b_2, \ldots, b_n)^T$. By Lemma 2.1, we get $A_n\left(\frac{1}{2}\right) \neq 0$. Then, by the Crammer rule, (2.2) has a unique solution for $a_0, a_1, a_2, \ldots, a_n$ in term of $b_0, b_1, b_2, \ldots, b_n$, which satisfies

$$|a_k| = \frac{\left|A_n^k\left(\frac{1}{2}\right)\right|}{\left|A_n\left(\frac{1}{2}\right)\right|} \le \frac{n!b}{\left|A_n\left(\frac{1}{2}\right)\right|}, \text{ for } k = 0, 1, 2, \ldots, n. \tag{2.3}$$

This implies that Coe($D$) is $|\cdot|$-bounded by $\frac{n!b}{\left|A_n\left(\frac{1}{2}\right)\right|}$, which does not depend on $p \in D$. □

**Corollary 2.3**. *Let $n$ be a positive integer. Let $D$ be a nonempty subset of $\mathcal{P}_n$ with coefficient set* Coe($D$). *Define $D' = \{p' \in \mathcal{P}_{n-1}: p \in D\}$. Then, we have*

$$D \textit{ is } \|\cdot\|\textit{-bounded} \quad \Leftrightarrow \quad D' \textit{ is } \|\cdot\|\textit{-bounded}.$$

*Proof*. Suppose that $D$ is $\|\cdot\|$-bounded by $b$, for some $b > 0$. By (2.3) in the proof of Proposition 2.2, Coe($D$) is $|\cdot|$-bounded by $\frac{n!b}{\left|A_n\left(\frac{1}{2}\right)\right|}$. This implies that Coe($D'$) is $|\cdot|$-bounded by $\frac{n!nb}{\left|A_n\left(\frac{1}{2}\right)\right|}$. By Proposition 2.2 again, $D'$ is $\|\cdot\|$-bounded by $\frac{n!n^2b}{\left|A_n\left(\frac{1}{2}\right)\right|}$. □

## 3. The Metric Projection $P_{\mathcal{P}_n}$: $C[0, 1] \to \mathcal{P}_n$ Is a Single-Valued Mapping

In this section, we recall some results from [16] regarding the metric projection operator $P_{\mathcal{P}_n}$. One of the important results obtained in [16] is that $P_{\mathcal{P}_n}$ is a single-valued mapping, which is proved by using Corollary 2.3 in the previous section and the Chebyshev's Equioscillation Theorem. We first recall some notations defined in [16]. For any $f \in C[0, 1]$, define

$$M(f) = \{t \in [0, 1]: |f(t)| = \|f\|\}.$$

By the continuity of $f$ on [0, 1], $M(f)$ is a nonempty closed subset of [0, 1], which is called the maximizing set of $f$. The following lemma shows the connection between the normalized duality mapping $J$: $C[0, 1] \to 2^{C^*[0,1]}\backslash\{\emptyset\}$ and the maximizing sets, which are studied in [16].

**Lemma 3.1** (see [16]). *For any $f \in C[0, 1]$ with $\|f\| > 0$, then*

(a) $\mu \in J(f) \Longrightarrow v(\mu, [0, 1]\backslash M(f)) = 0$;

(b) *Let $\{t_j \in M(f): j = 1, 2, \dots, m\} \subseteq M(f)$. For a positive integer m, define $\mu \in rca[0, 1]$ by*

$$\mu(t_j) = \alpha_j f(t_j), \text{ for } j = 1, 2, \dots, m \quad \textit{and} \quad \mu([0, 1]\backslash\{t_j: j = 1, 2, \dots, m\}) = 0,$$

*where $\alpha_j > 0$, for $j = 1, 2, \dots, m$ satisfying $\sum_{j=1}^{m} \alpha_j = 1$. Then, $\mu \in J(f)$.*

(c) *If $M(f)$ is not a singleton, then $J(f)$ is an infinite set.*

Let $f \in C[0, 1]$ and $p \in \mathcal{P}_n$. Then $p \in P_{\mathcal{P}_n}(f)$ if and only if $\|f - p\| = \min\{\|f - q\|: q \in \mathcal{P}_n\}$.

That is, if $p \in P_{\mathcal{P}_n}(f)$, then $p$ is the polynomial of best uniformly approximation for $f$ satisfying,

$$\begin{aligned}
&\|f - p\| \\
&= \max_{0 \le t \le 1} |f(t) - p(t)| \\
&= \min\{\|f - q\|: q(t) = \textstyle\sum_{k=0}^{n} c_k t^k \in \mathcal{P}_n\} \\
&= \min\{\max_{0 \le t \le 1} |f(t) - \textstyle\sum_{k=0}^{n} c_k t^k|: \sum_{k=0}^{n} c_k t^k \in \mathcal{P}_n\}.
\end{aligned}$$

The metric projection operator $P_{\mathcal{P}_n}$ from $C[0, 1]$ to $\mathcal{P}_n$ has the following existence properties, which have been proved in [16]. We will reprove it here with more details about the properties of bounded subset in $\mathcal{P}_n$, which are proved in section 2.

**Proposition 3.2 (see [16])**. *Let n be a nonnegative integer. Then*

(a) $P_{\mathcal{P}_n}(p) = p$, *for any* $p \in \mathcal{P}_n$;

(b) $P_{\mathcal{P}_n}(f) \neq \emptyset$, *for any* $f \in C[0, 1]$.

*Proof.* Part (a) is clear. We only prove part (b) for $f \in C[0, 1]\backslash\mathcal{P}_n$. Since $\mathcal{P}_n$ is a closed subspace of $C[0, 1]$, then, for any $f \in C[0, 1]\backslash\mathcal{P}_n$, we have

$$A \coloneqq \min\{\|f - q\|: q \in \mathcal{P}_n\} > 0.$$

We take a sequence $\{q_m\} \subseteq \mathcal{P}_n$ such that

$$\|f - q_m\| \downarrow A, \text{ as } m \to \infty.$$

This implies that $\{\|q_m\|\}$ is bounded. Since every $q'_m$ is a polynomial with degree less than or equal to $n - 1$, by the boundness of $\{\|q_m\|\}$ and by Corollary 2.3, it follows that $\{\|q'_m\|\}$ is bounded. We define the upper bound $d = \max\{\|q'_m\|: m = 1, 2, 3 \ldots\}$. For any given $\varepsilon > 0$, we take $\delta > 0$ with $\delta < \frac{\varepsilon}{d}$. By the mean value theorem, this implies that, for any $s, t \in [0, 1]$, we have

$$|s - t| < \delta \implies |q_m(s) - q_m(t)| < \varepsilon, \text{ for } m = 1, 2, \ldots .$$

This implies that $\{\|q_m\|\}$ is a bounded and equicontinuous subset in $C[0, 1]$. By Azsela-Ascali Theorem, $\{\|q_m\|\}$ is conditionally compact. Hence, there is a subsequence $\{q_{m_i}\} \subseteq \{q_m\} \subseteq \mathcal{P}_n$ and there is $q \in \mathcal{P}_n$ such that $q_{m_i} \to q$, as $i \to \infty$. We have

$$\|f - q\| \leq \|f - q_{m_i}\| + \|q - q_{m_i}\| \to A, \text{ as } i \to \infty.$$

This implies $\|f - q\| = A$. That is $q \in P_{\mathcal{P}_n}(f)$. □

In 1859, Chebyshev proved the well-known theorem called the Chebyshev's Equioscillation theorem. This theorem provides the criterion for the necessary and sufficient conditions for the solutions of best approximation of continuous functions by polynomials with degree less than or equal to $n$. We list the revised version of the Chebyshev's Equioscillation Theorem with respect to the metric projection $P_{\mathcal{P}_n}$ (see [5, 6, 29, 30]).

**Chebyshev's Equioscillation Theorem**. *Let* $f(t)$ *be a continuous function on* $[0,1]$ *and, for* $p(t) = \sum_{k=0}^{n} a_k t^k \in \mathcal{P}_n$, *suppose that*

$$A(f, p) \coloneqq \|f - p\| = \max_{0 \leq t \leq 1} |f(t) - p(t)|.$$

*Then,* $p \in P_{\mathcal{P}_n}(f)$ *if and only if, there are* $n+2$ *points* $0 \leq t_0 < t_1 < \cdots < t_{n+1} \leq 1$ *such that*

$$f(t_i) - p(t_i) = \epsilon A(f, p)(-1)^i, \; i = 0, 1, 2, \ldots, n+1, \tag{2.3}$$

*where* $\epsilon = 1$ *or* $-1$.

The set $\{t_0, t_1, \ldots, t_{n+1}\}$ satisfying (2.3) is called a $n$-Chebyshev set of $f$ with respect to $p$. The collection of all $n$-Chebyshev set of $f$ with respect to $p$, is denoted by $S(f, p)$. It is clear that

$$S(f, p) \subseteq M(f - p), \text{ for } p \in P_{\mathcal{P}_n}(f).$$

In the following example, we show that $S(f, p)$ may not be unique, in general.

**Example 3.3** (see [16]). Let $n = 1$. Let $f(t) = \sin(4\pi t) \in C[0, 1]$. Then, $P_{\mathcal{P}_1}(f) = \theta$ with $A(f, \theta) = 1$. One has that there are two 1-Chebyshev sets of $f$ with respect to $p$, which are

$$\left\{\frac{\pi}{2}, \frac{3\pi}{2}, \frac{5\pi}{2}\right\} \text{ with } \epsilon = 1 \quad \text{and} \quad \left\{\frac{3\pi}{2}, \frac{5\pi}{2}, \frac{7\pi}{2}\right\} \text{ with } \epsilon = -1.$$

Since the Banach space $C[0, 1]$ is neither reflexive, nor strictly convex and smooth, so, it is not for sure that the metric projection $P_{\mathcal{P}_n}$ from $C[0, 1]$ to $\mathcal{P}_n$ is a single-valued mapping or not. However, by Proposition 2.2, the following theorem shows that $\mathcal{P}_n$ is indeed a single-valued mapping (See [5]). However, [16] provides a directly proof of the following theorem.

**Theorem 3.4 (see [16])**. *Let n be a nonnegative integer. The metric projection* $P_{\mathcal{P}_n}: C[0, 1] \to \mathcal{P}_n$ *is a single-valued mapping.*

By the Chebyshev's Equioscillation Theorem, we have the following approximation properties.

**Proposition 3.5** (see [16]. *Let n be a nonnegative integer. Let* $f \in C[0, 1]$ *with* $p = P_{\mathcal{P}_n}(f)$. *Then,*

(i) $\beta p = P_{\mathcal{P}_n}(\beta f)$, *for any* $\beta \in \mathbb{R}$,

$A(\beta f, \beta p) = |\beta| A(f, p)$ *and* $S(\beta f, \beta p) = S(f, p)$,

$\epsilon(\beta f, \beta p) = \epsilon(f, p)$, *if* $\beta \geq 0$ *and* $\epsilon(\beta f, \beta p) = -\epsilon(f, p)$, *if* $\beta < 0$.

(ii) $P_{\mathcal{P}_n}(f + q) = P_{\mathcal{P}_n}(f) + q = p + q$, *for any* $q \in \mathcal{P}_n$,

$A(f + q, p + q) = A(f, p)$ *and* $S(f + q, p + q) = S(f, p)$.

*In particular, we have*

$$p + \lambda = P_{\mathcal{P}_n}(f + \lambda), \text{ for any } \lambda \in \mathbb{R};$$

(iii) *Let* $f \in C[0, 1] \backslash \mathcal{P}_n$ *and* $p \in \mathcal{P}_n$ *with* $p = P_{\mathcal{P}_n}(f)$. *Then*

$$p = P_{\mathcal{P}_n}(\alpha f + (1 - \alpha)p), \text{ for any } \alpha \in (0, 1).$$

## 4. The Continuity of Metric Projection Operator $P_{\mathcal{P}_n}: C[0, 1] \to \mathcal{P}_n$

In the previous section, we proved that the metric projection operator $P_{\mathcal{P}_n}: C[0, 1] \to \mathcal{P}_n$ is a single-valued mapping. One step further, in this section, we prove that $P_{\mathcal{P}_n}$ is also (norm to norm) continuous, which is the main theorem of this paper.

**Theorem 4.1**. *Let n be a positive integer. The metric projection* $P_{\mathcal{P}_n}: C[0, 1] \to \mathcal{P}_n$ *is a single-valued continuous mapping.*

*Proof.* By Theorem 3.4, $P_{\mathcal{P}_n}: C[0, 1] \to \mathcal{P}_n$ is a single-valued mapping. Let $\{f_m\}$ be a convergent sequence in $C[0, 1]$ with limit $f$. Then $\{f_m\}$ is $\|\cdot\|$-bounded. Let $P_{\mathcal{P}_n}(f_m) = p_m$, for all $m$. Take an arbitrary $w \in \mathcal{P}_n$, we have

$$\|p_m\| \leq \|p_m - f_m\| + \|f_m\| \leq \|w - f_m\| + \|f_m\| \leq \|w\| + 2\|f_m\|.$$

This implies that $\{p_m\}$ is $\|\cdot\|$-bounded. By Corollary 2.3, $\{p'_m\}$ is $\|\cdot\|$-bounded. We define the upper bound $d = \max\{\|p'_m\|: m = 1, 2, 3 \ldots\}$. For any given $\varepsilon > 0$, we take $\delta > 0$ with $\delta < \frac{\varepsilon}{d}$. By the mean value theorem, this implies that, for any $s, t \in [0, 1]$, we have

$$|s-t| < \delta \implies |p_m(s) - p_m(t)| < \varepsilon, \text{ for every } m = 1, 2, \ldots .$$

This implies that $\{p_m\}$ is a bounded and equicontinuous subset in $C[0, 1]$. By Azsela-Ascali Theorem, $\{p_m\}$ is conditionally compact (see [7]). Since $\mathcal{P}_n$ is closed in $C[0, 1]$, then, there is a subsequence $\{p_{m_i}\} \subseteq \{p_m\} \subseteq \mathcal{P}_n$ and there is $p \in \mathcal{P}_n$ such that $p_{m_i} \to p$, as $i \to \infty$. For any $q \in \mathcal{P}_n$, by $P_{\mathcal{P}_n}(f_m) = p_m$, we have

$$\begin{aligned} \|f-p\| &\le \|f-f_{m_i}\| + \|f_{m_i} - p_{m_i}\| + \|p - p_{m_i}\| \\ &\le \|f-f_{m_i}\| + \|f_{m_i} - q\| + \|p - p_{m_i}\| \\ &\le \|f-f_{m_i}\| + \|f_{m_i} - f\| + \|f - q\| + \|p - p_{m_i}\|. \end{aligned}$$

By $p_{m_i} \to p$ and $f_{m_i} \to f$ as $i \to \infty$, this implies

$$\|f-p\| \le \|f-q\|, \text{ for any } q \in \mathcal{P}_n.$$

It follows that $P_{\mathcal{P}_n}(f) = p$. This proves the continuity of $P_{\mathcal{P}_n}$. □

## 5. Mordukhovich Derivatives and Gâteaux Directional Derivatives of $P_{\mathcal{P}_n}$ in $C[0, 1]$

In the analysis of Banach spaces, Gâteaux directional differentiation and Fréchet differentiation are only deal with single-valued mappings in Banach spaces. In particular, the differentiability of the metric projection operator has been studied by many authors in Hilbert spaces or in uniformly convex and uniformly smooth Banach spaces (see [3, 4, 8, 10, 11, 13−18, 22, 26, 27]).

It should be congratulated on the theory of generalized differentiation, which extends the Fréchet differentiation of single-valued mappings to Mordukhovich differentiation of set-valued mappings. This new theory has been rapidly developed and has been widely applied to several branches of mathematics such as, operator theory, optimization theory, approximation theory, control theory, equilibrium theory, and so forth (see [19−21]).

In this section, we use the continuity of the metric projection operator $P_{\mathcal{P}_n}$ to study some properties of Mordukhovich derivatives of $P_{\mathcal{P}_n}$. We first recall some concepts related to Mordukhovich derivatives.

Let $(X, \|\cdot\|_X)$ and $(Y, \|\cdot\|_Y)$ be real Banach spaces with topological dual spaces $X^*$ and $Y^*$, respectively. Let $\langle\cdot, \cdot\rangle_X$ denote the real canonical pairing between $X^*$ and $X$ and $\langle\cdot, \cdot\rangle_Y$ the real canonical pairing between $Y^*$ and $Y$. Let $\Delta$ be a nonempty subset of $X$ and let $F\colon \Delta \rightrightarrows Y$ be a set valued mapping. The graph of $F$ is defined by the following subset in $\Delta \times Y$

$$\text{gph}F = \{(x, y) \in \Delta \times Y : y \in F(x)\}.$$

For $(x, y) \in \text{gph}F$, that is, for $x \in \Delta$ and $y \in F(x)$, the Mordukhovich derivative (which is also called Mordukhovich coderivative, or coderivative) of $F$ at point $(x, y)$ is a set valued mapping $\widehat{D}^*F(x,y): Y^* \rightrightarrows X^*$. For any $y^* \in Y^*$, it is defined by (see Definitions 1.13 and 1.32 in Chapter 1 in [19])

$$\widehat{D}^*F(x,y)(y^*) = \left\{ z^* \in X^*: \limsup_{\substack{(u,v)\to(x,y) \\ (u,v)\in \text{gph}F}} \frac{\langle z^*, u-x\rangle_X - \langle y^*, v-y\rangle_Y}{\|u-x\|_X + \|v-y\|_Y} \leq 0 \right\}. \tag{5.1}$$

If $(x,y) \notin \text{gph}F$, then, we define

$$\widehat{D}^*F(x,y)(y^*) = \emptyset, \text{ for any } y^* \in Y^*.$$

By the above definition (5.1), $\widehat{D}^*F(x,y): Y^* \rightrightarrows X^*$ is a set valued mapping, which is called the Mordukhovich differential operator (or the Mordukhovich codifferential operator) of $F$ at $(x,y)$.

In particular, let $g: \Delta \to Y$ be a single-valued continuous mapping. By (5.1), the Mordukhovich derivative of $g$ at point $(x, g(x))$ is a set valued mapping $\widehat{D}^*g(x,g(x)): Y^* \rightrightarrows X^*$. In this paper, $\widehat{D}^*g(x,g(x))$ is simply written as $\widehat{D}^*g(x)(y^*)$, for single-valued continuous mapping $g$. For any $x \in \Delta$ and $y^* \in Y^*$, it is defined by

$$\widehat{D}^*g(x,g(x))(y^*) = \widehat{D}^*g(x)(y^*) = \left\{ z^* \in X^*: \limsup_{\substack{u\to x \\ u\in\Delta}} \frac{\langle z^*, u-x\rangle_X - \langle y^*, g(u)-g(x)\rangle_Y}{\|u-x\|_X + \|g(u)-g(x)\|_Y} \leq 0 \right\}. \tag{5.2}$$

By Theorem 4.1, the metric projection operator $P_{\mathcal{P}_n}: C[0,1] \to \mathcal{P}_n$ is a single-valued continuous mapping. From the definition in (5.2), in particular, for any $f \in C[0,1]$ with $p = P_{\mathcal{P}_n}(f)$ and for $\mu \in C^*[0,1]$, we have

$$\widehat{D}^*P_{\mathcal{P}_n}(f)(\mu) = \left\{ \varphi \in C^*[0,1]: \limsup_{\substack{g\to f \\ g\in C[0,1]}} \frac{\langle \varphi, g-f\rangle - \langle \mu, P_{\mathcal{P}_n}(g) - P_{\mathcal{P}_n}(f)\rangle}{\|g-f\| + \|P_{\mathcal{P}_n}(g) - P_{\mathcal{P}_n}(f)\|} \leq 0 \right\}$$

$$= \left\{ \varphi \in C^*[0,1]: \limsup_{\substack{g\to f \\ g\in C[0,1]}} \frac{\langle \varphi, g-f\rangle - \langle \mu, P_{\mathcal{P}_n}(g) - p\rangle}{\|g-f\| + \|P_{\mathcal{P}_n}(g) - p\|} \leq 0 \right\}. \tag{5.3}$$

**Theorem 5.1**. *Let $n$ be a positive integer. Let $f \in C[0,1]$ and $p \in \mathcal{P}_n$ with $p = P_{\mathcal{P}_n}(f)$. For any $\mu, \gamma \in C^*[0,1]$, we have*

(i) $\mu([0,1]) \neq \gamma([0,1]) \implies \gamma \notin \widehat{D}^*P_{\mathcal{P}_n}(f)(\mu)$ *and* $\mu \notin \widehat{D}^*P_{\mathcal{P}_n}(f)(\gamma)$;

(ii) $\mu([0,1]) \neq 0 \implies \theta^* \notin \widehat{D}^*P_{\mathcal{P}_n}(f)(\mu)$ *and* $\mu \notin \widehat{D}^*P_{\mathcal{P}_n}(f)(\theta^*)$;

(iii) $\langle \gamma, f - p\rangle < 0 \implies \gamma \notin \widehat{D}^*P_{\mathcal{P}_n}(f)(\mu)$, *for any* $\mu \in C^*$.

*Proof*. Proof of (i). For any real number $\lambda$, we consider it as a constant function defined on [0, 1] with value $\lambda$. Then, for any $\mu \in C^*$, we have

$$\langle\mu, \lambda\rangle = \lambda\mu([0,1]), \text{ for any } \lambda \in \mathbb{R}. \quad (5.4)$$

For any $\mu, \gamma \in C^*$, by definition, we have

$$\gamma \in \widehat{D}^* P_{\mathcal{P}_n}(f)(\mu) \Leftrightarrow \limsup_{\substack{(g,q)\to(f,p) \\ q=P_{\mathcal{P}_n}(g)}} \frac{\langle\gamma, g-f\rangle - \langle\mu, q-p\rangle}{\|g-f\| + \|q-p\|} \le 0. \quad (5.5)$$

Suppose that $\mu([0,1]) < \gamma([0,1])$. Then, in limit (5.5), we take a line segment direction as $g_\lambda = f + \lambda$, for $\lambda \downarrow 0$. One has

$$g_\lambda - f = \lambda \to \theta, \text{ as } \lambda \downarrow 0.$$

By **Proposition 3.5**, we have

$$p + \lambda = P_{\mathcal{P}_n}(f + \lambda) = P_{\mathcal{P}_n}(g_\lambda) \coloneqq q_\lambda, \text{ for any } \lambda > 0.$$

It is clear that $p + \lambda \to p$, as $\lambda \downarrow 0$. By (5.4), we calculate the limit in (5.5),

$$\limsup_{\substack{(g,q)\to(f,p) \\ q=P_{\mathcal{P}_n}(g)}} \frac{\langle\gamma, g-f\rangle - \langle\mu, q-p\rangle}{\|g-f\| + \|q-p\|}$$

$$\ge \limsup_{\substack{(g_\lambda,q_\lambda)\to(f,p) \\ q_\lambda=P_{\mathcal{P}_n}(g_\lambda)}} \frac{\langle\gamma, g_\lambda-f\rangle - \langle\mu, q_\lambda-p\rangle}{\|g_\lambda-f\| + \|q_\lambda-p\|}$$

$$= \limsup_{(f+\lambda,p+\lambda)\to(f,p)} \frac{\langle\gamma, f+\lambda-f\rangle - \langle\mu, p+\lambda-p\rangle}{\|f+\lambda-f\| + \|p+\lambda-p\|}$$

$$= \limsup_{\lambda\downarrow 0} \frac{\langle\gamma,\lambda\rangle - \langle\mu,\lambda\rangle}{2\lambda}$$

$$= \limsup_{\lambda\downarrow 0} \frac{\lambda(\gamma([0,1]-\mu([0,1])}{2\lambda}$$

$$= \frac{\gamma([0,1]-\mu([0,1]}{2} > 0.$$

By (5.5), this concludes that

$$\mu([0,1]) < \gamma([0,1]) \implies \gamma \notin \widehat{D}^* P_{\mathcal{P}_n}(f,p)(\mu). \quad (5.6)$$

Similarly, to the proof of (5.6), in the limit (5.5), by taking $g_\lambda = f - \lambda$, for $\lambda \downarrow 0$, we can show

$$\mu([0,1]) > \gamma([0,1]) \implies \gamma \notin \widehat{D}^* P_{\mathcal{P}_n}(f,p)(\mu).$$

This proves part (i) of this theorem. Part (ii) is proved by following from part (i) immediately.

Proof of (iii). Suppose $\langle\gamma, f - p\rangle < 0$. In the limit (5.5), we take a line segment direction as $g_\alpha = (1-\alpha)f + \alpha p$, for $\alpha \downarrow 0$ with $\alpha < 1$. By part (iii) in Proposition 3.5, one has

$$p = P_{\mathcal{P}_n}((1-\alpha)f + \alpha p), \text{ for any } \alpha \in (0, 1).$$

Since $(1-\alpha)f + \alpha p \to f$, as $\alpha \downarrow 0$ with $\alpha < 1$, we have

$$\limsup_{\substack{(g,q)\to(f,p) \\ q=P_{\mathcal{P}_n}(g)}} \frac{\langle \gamma, g-f\rangle - \langle \mu, q-p\rangle}{\|g-f\| + \|q-p\|}$$

$$\geq \limsup_{\substack{(g_\alpha,q_\alpha)\to(f,p) \\ q_\alpha=P_{\mathcal{P}_n}(g_\alpha)}} \frac{\langle \gamma, g_\alpha-f\rangle - \langle \mu, q_\alpha-p\rangle}{\|g_\alpha-f\| + \|q_\alpha-p\|}$$

$$= \limsup_{\alpha\downarrow 0, \alpha<1} \frac{\langle \gamma, (1-\alpha)f + \alpha p - f\rangle - \langle \mu, p-p\rangle}{\|(1-\alpha)f + \alpha p - f\| + \|p-p\|}$$

$$= \limsup_{\alpha\downarrow 0, \alpha<1} \frac{-\langle \gamma, \alpha(f-p)\rangle}{\|\alpha(f-p)\|}$$

$$= \frac{-\langle \gamma, f-p\rangle}{\|f-p\|} > 0.$$

By (5.5), from this, this concludes that

$$\langle \gamma, f-p\rangle < 0 \implies \gamma \notin \widehat{D}^* P_{\mathcal{P}_n}(f,p)(\mu), \text{ for any } \mu \in C^*. \qquad \square$$

In [12], some fixed-point properties of the Mordukhovich differential operator are proved in uniformly convex and uniformly smooth Banach spaces. In this section, we study the properties of fixed-point properties of the Mordukhovich differential operator in the Banach space $C[0, 1]$, which is not uniformly convex and uniformly smooth.

Recall that, for any $f \in C[0, 1]$ and $p \in \mathcal{P}_n$ with $p = P_{\mathcal{P}_n}(f)$, $\widehat{D}^* P_{\mathcal{P}_n}(f): C^*[0,1] \rightrightarrows C^*[0,1]$ is a set-valued mapping. Let $\mathcal{F}\left(\widehat{D}^* P_{\mathcal{P}_n}(f)\right)$ denote the collection of all fixed points of Mordukhovich differential operator $\widehat{D}^* P_{\mathcal{P}_n}(f)$. Then, by (5.3), we have

$$\mathcal{F}\left(\widehat{D}^* P_{\mathcal{P}_n}(f)\right) = \left\{\varphi \in C^*[0,1]: \limsup_{g\to f} \frac{\langle \varphi, (g-f)\rangle - \langle \varphi, P_{\mathcal{P}_n}(g) - p\rangle}{\|g-f\| + \|P_{\mathcal{P}_n}(g) - p\|} \leq 0\right\}$$

$$= \left\{\varphi \in C^*[0,1]: \limsup_{g\to f} \frac{\langle \varphi, (g-f) - (P_{\mathcal{P}_n}(g) - p)\rangle}{\|g-f\| + \|P_{\mathcal{P}_n}(g) - p\|} \leq 0\right\}.$$

We write the orthogonal space of $\mathcal{P}_n$ by $\mathcal{P}_n^\perp$, which is defined as follows.

$$\mathcal{P}_n^\perp = \{\mu \in C^*[0,1]: \langle \mu, p\rangle = 0, \text{for every } p \in \mathcal{P}_n\}.$$

**Theorem 5.2**. *Let $n$ be a positive integer. Let $f \in C[0, 1]$ with $p = P_{\mathcal{P}_n}(f)$. Let $\mu \in C^*[0, 1]$ with $\langle \mu, f\rangle \neq 0$, we have*

(i) $\gamma \in \mathcal{P}_n^\perp \implies \mu \notin \widehat{D}^* P_{\mathcal{P}_n}(f)(\gamma)$, *for* $\gamma \in C^*[0, 1]$;

(ii) $\mu \in \mathcal{P}_n^\perp \implies \mu \notin \mathcal{F}\left(\widehat{D}^* P_{\mathcal{P}_n}(f)\right)$.

*Proof*. Proof of (i). For the given $f \in C[0, 1]$ and $p = P_{\mathcal{P}_n}(f)$, let $\mu \in C^*[0, 1]$ with $\langle \mu, f\rangle \neq 0$. Without loose of the generality, we assume that $\langle \mu, f\rangle > 0$. Let $\gamma \in \mathcal{P}_n^\perp$. We take a line segment direction as $g_\lambda = f + \lambda f$, for $\lambda \downarrow 0$. One has

$$g_\lambda = f + \lambda f \to f, \text{ as } \lambda \downarrow 0.$$

By Proposition 3.5, we have that $P_{\mathcal{P}_n}(g_\lambda) = P_{\mathcal{P}_n}(f + \lambda f) = (1+\lambda)\, p$. This implies that $P_{\mathcal{P}_n}(g_\lambda) \to p$, as $\lambda \downarrow 0$. By $\gamma \in \mathcal{P}_n^{\perp}$, we calculate

$$\limsup_{\substack{(g,q)\to(f,p)\\ q=P_{\mathcal{P}_n}(g)}} \frac{\langle \mu, g-f\rangle - \langle \gamma,\ q-p\rangle}{\|g-f\| + \|q-p\|}$$

$$\geq \limsup_{(g_\lambda, P_{\mathcal{P}_n}(g_\lambda))\to(f,p)} \frac{\langle \mu, g_\lambda - f\rangle - \langle \gamma,\ P_{\mathcal{P}_n}(g_\lambda) - p\rangle}{\|g_\lambda - f\| + \|P_{\mathcal{P}_n}(g_\lambda) - p\|}$$

$$= \limsup_{(f+\lambda h, P_{\mathcal{P}_n}(f+\lambda h))\to(f,p)} \frac{\lambda\langle \mu, f\rangle}{\lambda\|f\| + \lambda\|p\|}$$

$$= \frac{\langle \mu, f\rangle}{\|f\| + \|p\|} > 0.$$

This implies that $\mu \notin \widehat{D}^* P_{\mathcal{P}_n}(f, p)(\gamma)$. Then, part (ii) follows from part (i) immediately. □

By Theorem 4.1 in the previous section, the metric projection operator $P_{\mathcal{P}_n}: C[0, 1] \to \mathcal{P}_n$ is a single valued continuous mapping. So, this allows us to consider the Gâteaux directional differentiability of the metric projection operator $P_{\mathcal{P}_n}$.

**Proposition 5.3.** *Let n be a nonnegative integer. We have*

(i) *For any* $f \in C[0, 1]$,

$$P'_{\mathcal{P}_n}(f)(q) = q, \text{ for any } q \in \mathcal{P}_n \backslash \{\theta\};$$

(ii) *For any* $q \in \mathcal{P}_n$,

$$P'_{\mathcal{P}_n}(q)(f) = P_{\mathcal{P}_n}(f), \text{ for any } f \in C[0, 1] \backslash \{\theta\}.$$

*Proof*. Proof of (i). For any $f \in C[0, 1]$ and for any $q \in \mathcal{P}_n \backslash \{\theta\}$, by the properties of the metric projection operator $P_{\mathcal{P}_n}$ in Proposition 3.5 and by $tq \in \mathcal{P}_n \backslash \{\theta\}$, for any $t > 0$, we have

$$P'_{\mathcal{P}_n}(f)(q) = \lim_{t\downarrow 0} \frac{P_{\mathcal{P}_n}(f+tq) - P_{\mathcal{P}_n}(f)}{t}$$

$$= \lim_{t\downarrow 0} \frac{P_{\mathcal{P}_n}(f) + tq - P_{\mathcal{P}_n}(f)}{t} = q.$$

Proof of (ii). For any $q \in \mathcal{P}_n$ and $f \in C[0, 1] \backslash \{\theta\}$, by the properties of the metric projection operator $P_{\mathcal{P}_n}$ in Proposition 3.5 and by $P_{\mathcal{P}_n}(tf) = tP_{\mathcal{P}_n}(f) \in \mathcal{P}_n$, for any $t > 0$, we have

$$P'_{\mathcal{P}_n}(q)(f) = \lim_{t\downarrow 0} \frac{P_{\mathcal{P}_n}(q+tf) - P_{\mathcal{P}_n}(q)}{t}$$

$$= \lim_{t\downarrow 0} \frac{q + P_{\mathcal{P}_n}(tf) - q}{t}$$

$$= \lim_{t\downarrow 0} \frac{tP_{\mathcal{P}_n}(f)}{t}$$

$$= P_{\mathcal{P}_n}(f). \qquad \square$$